\documentclass[a4paper]{amsart}
\usepackage{amssymb,bbm,amsmath,amsfonts,amsthm, amscd}%,amsbsy}
\usepackage{xcolor}%,refcheck}
\usepackage{lineno}
\usepackage[linktocpage,colorlinks=true,raiselinks]{hyperref}
\hypersetup{urlcolor=blue, citecolor=red, linkcolor=blue}

\newtheorem{definition}{Definition}[section]

\newtheorem{lemma}{Lemma}[section]
\newtheorem{theorem}{Theorem}[section]

\newtheorem{remark}{Remark}[section]
\newcommand{\br}{r}
\newcommand{\bu}{u}
\newcommand{\bv}{v}
\newcommand{\bvn}{v}
\newcommand{\bvnh}{\widetilde{v}}
\newcommand{\bve}{v^E}
\newcommand{\bvee}{\bve_{\vep}}
\newcommand{\bx}{x}
\newcommand{\by}{y}

\newcommand{\bphi}{\varphi}
\newcommand{\vep}{\varepsilon}
\newcommand{\bC}{C}
\newcommand{\de}{\mathrm{d}}
\newcommand{\dx}{\,\de\bx}
\newcommand{\dy}{\,\de\by}

\newcommand{\N}{\mathbbm{N}}

\newcommand{\R}{\mathbbm{R}}
\newcommand{\T}{\mathbbm{T}}

\newcommand{\tore}{{\mathbbm{T}^{3}}}
\DeclareMathOperator{\dive}{div}
\begin{document}

%\begin{frontmatter}

  \title[Energy conservation for 
    incompressible fluids: the H\"older case]{Energy conservation for weak solutions of
    incompressible fluid equations: the H\"older case and connections with
    Onsager's conjecture} 
  \author{Luigi  C. Berselli}
  \address{Dipartimento di Matematica, Universit\`a di Pisa, Via
    F. Buonarroti 1/c, Pisa, I56127, ITALY}
\email{luigi.carlo.berselli@unipi.it}
  \begin{abstract}
    In this paper we give elementary proofs of energy conservation for
    weak solutions to the Euler and Navier-Stokes equations in the
    class of H\"older continuous functions, relaxing some of the
    assumptions on the time variable (both integrability and
    regularity at initial time) and presenting them in a unified way. 

    Then, in the final section we prove (for the Navier-Stokes equations) a result of
    energy conservation in presence of a solid boundary and with Dirichlet boundary
    conditions. This result seems the first one --in the viscous case-- with H\"older type
    assumptions, but without additional assumptions on the pressure.
%\end{keyword}
  \end{abstract}
  
%\begin{keyword}
%Energy conservation, Onsager's conjecture, Euler and Navier-Stokes equations
%\MSC[2010] 35Q30
%\end{keyword}

 \maketitle
% \end{frontmatter}

%\linenumbers

\section{Introduction}
In this paper we start proving (and we slightly improve) in an elementary way
--that is without resorting to Paley-Littlewood decomposition or other
tools from harmonic analysis-- results about energy conservation for
weak solutions to incompressible fluids.  We consider weak solutions
to the Euler equations
\begin{equation}
  \label{eq:Ebvip}
  \begin{aligned}
    \partial_{t}\bve+(\bve\cdot\nabla)\,\bve+\nabla p^{E}&=0\qquad
    &(t,\bx)\in (0,T)\times \tore,
    \\
    \dive\bve&=0\qquad &(t,\bx)\in (0,T)\times \tore,
    \\
    \bve(0,\bx)&=\bve_0(\bx) \qquad &\bx\in \tore,
  \end{aligned}
\end{equation}
or Leray-Hopf weak solutions of the Navier-Stokes equations (NSE in
the sequel)
\begin{equation}
  \label{eq:NSEbvip}
  \begin{aligned}
    \partial_{t}\bvn+(\bvn\cdot\nabla)\,\bvn-\nu\Delta\bvn+\nabla
    p&=0\qquad &(t,\bx)\in (0,T)\times \tore,
    \\
    \dive\bvn&=0\qquad &(t,\bx)\in (0,T)\times \tore,
    \\
    \bvn(0,\bx)&=\bvn_0(\bx) \qquad &\bx\in \tore,
  \end{aligned}
\end{equation}
for fixed $\nu>0$. We begin from the space-periodic case, under the
additional assumption of having a $C^{0,\alpha}(\tore)$ velocity, as
suggested by Onsager conjecture~\cite{Ons1949}.  Note that the
boundary value problem for the Euler equations can be treated at the
price of some technical steps, but more or less in the same way, by
adapting results from~\cite{BT2018}. In the final section we will also
prove some results in the case of a domain with solid boundaries
and Dirichlet boundary conditions for the NSE for
which the treatment of the pressure is problematic.  In the
sequel we will us assume $f=0$ for simplicity, but results can be
easily adapted to consider also non-zero external forces.

The Onsager conjecture (only recently solved) suggested the value
$\alpha=1/3$, especially for the case of the Euler
equations~\eqref{eq:Ebvip}, but the conjecture was mainly considering
only the H\"older regularity with respect the \textit{space
  variables}. Here, we consider a combination of space-time
conditions, proving families of criteria depending on the H\"older
exponent.

The first rigorous results about Onsager conjecture are probably those
of Eyink~\cite{Eyi1994,Eyi1995} and Constantin, E, Titi~\cite{CET1994}
in Fourier setting and in Besov spaces (slightly larger than H\"older
ones), respectively. A well known result is that if $\bve$ is a weak
solution to \eqref{eq:Ebvip} is such that
\begin{equation*}
   \bve\in L^{3}(0,T;B^{\alpha,\infty}_{3}(\tore))\cap
   C(0,T;L^{2}(\tore)),\qquad \text{with }\alpha>\frac{1}{3},
 \end{equation*}
 then $\|\bve(t)\|=\|\bve_{0}\|$, for all $t\in [0,T]$. As explained
 in \cite{CET1994} ``...This is basically the content of Onsager's
 conjecture, except Onsager stated his conjecture in H\"older spaces
 rather than Besov spaces. Obviously the above theorem implies similar
 results in Holder spaces...``
 \begin{remark}
   In this paper we want to consider the H\"older case, to keep the results the most
   elementary as possible, to be understood also by an audience familiar with classical
   spaces of mathematical analysis. This will also create links with the most recent
   developments of the theory, especially after the work of De Lellis and
   Sz{\'e}kelyhidi~\cite{DLS2009}, which originated an intense activity to prove the
   Onsager conjecture, mainly the non-conservation below the critical space $C^{1/3}$,
   with endpoint reached in~\cite{BDLSV2019,Ise2018}. We wish also to mention the recent
   result of Cheskidov and Luo~\cite{CL2020b}, again in the setting of H\"older continuous
   solutions, in connection with non-uniqueness, which treats scaling very similar to ones
   we obtain here. We avoid using Besov spaces (which can produce some technical
   improvements), both to keep the results elementary and to have a functional
   setting well-adapted also to the problem in a domain with boundary.
 \end{remark}
 The first five sections contain both a summary of known results and the
 simplification/extension of several theorems.  The most original part of the paper is
 presented in the last section, where we consider the Dirichlet problem for the NSE in the
 half-space. We prove energy conservation for a class of (partial) H\"older solutions,
 which is nevertheless larger than the scaling invariant class which provides uniqueness
 of weak solution and regularity of Leray solutions.

\vspace{.5cm}

 The results we prove here in the H\"older case follow easily from the
 method used in~\cite{CET1994}, and we are here trying to focus on the
 hypotheses in the time variable (instead of the spatial ones),
 showing how a slight more stringent assumption on the space variables
 (H\"older instead of Besov), produces some slight improvement in the
 time variable. We are not sure that all the results we prove in the
 first five sections could be considered as original (even if we could
 not find in the published literature), but nevertheless they are
 collected here in a unified way, providing also a summary of the
 state of the art.

 The result in~\cite{CET1994}, but in the setting of H\"older functions
 has been recently extended also to the boundary value problems for
 the Euler equations in Bardos and Titi~\cite{BT2018}, proving energy
 conservation for weak solutions such that
 \begin{equation*}
  \bve\in
  L^{3}(0,T;C^{0,\alpha}(\overline{\Omega})), \qquad \text{with
  }\alpha>\frac{1}{3},  
\end{equation*}
while technical improvements in the space variables can be found in
Beir\~ao da Veiga and Yiang~\cite{BY2020}. See also Duchon and
Robert~\cite{DR2000}.

We observe that the situation for the NSE is a little different and recent results are
those by Cheskidov~\cite{CFS2010} in the setting of standard (fractional) Sobolev spaces
of and Cheskidov Luo~\cite{CL2020}, see also Farwig and Taniuchi~\cite{FT2010}.  In the
case of Leray-Hopf weak solutions to NSE, the properties in the time variable are a little
improved and also the range of exponents is wider, if compared with the Euler case. The
sharpest known result in the Besov setting seems to the be the following one: Suppose that
$1\leq \beta<p\leq\infty$ are such that $\frac{2}{p}+\frac{1}{\beta}<1$. If $\bv$ is a
weak solution such that
\begin{equation}
    \label{eq:Cheskidov-Luo}
    \bvn\in
    L^{\beta}_{w}(0,T;B^{\frac{2}{p}+\frac{2}{\beta}-1}_{p,\infty}(\tore)),    
  \end{equation}
  then $\bv$ satisfies the energy equality.  This result contains as a
  particular case the classical ones by
  Prodi/Lions~\cite{Pro1959,Lio1960} and
  $\bvn\in L^{\beta}_{w}(0,T;B^{1/3}_{3,\infty}(\tore))$, for all
  $\beta>3$ is the borderline case. Recent extensions are those based
  on the condition on the gradient of $\bv$ in \cite{BC2020} and in
  Beir\~ao da Veiga and Yang~\cite{BY2019}. The results therein can be
  ``measured'' in terms of H\"older spaces as follows
  \begin{equation}
    \label{eq:BC-CJ}
    \bv \in L^{\frac{5}{3+2\alpha}}(0,T;C^{0,\alpha}(\overline{\Omega})),   
  \end{equation}
  which are the same as in~\eqref{eq:Cheskidov-Luo}in terms of
  scaling.  The results are in fact derived from
  $\nabla \bv\in L^{1+2/q}(0,T;L^{q}(\Omega))$, with $q>3$, by Morrey
  inclusion, which almost implies for instance
  $L^{\frac{15}{11}}(0,T;C^{0,1/3}(\overline{\Omega}))$. Nevertheless,
  results obtained by embedding --even if valid also for the boundary
  value problem-- are far away from being sharp. Hence, this is why we
  prove here results directly in the class of H\"older continuous
  functions.

  A traditional way to take advantage from the additional H\"older
  continuity of the solution is to use properties of mollification
  operators and not  Hilbert-space
  methods. For divergence-free functions, this is particularly
  delicate to be handled in the presence of a boundary.

  \bigskip

  \textbf{Plan of the paper:} In Section~\ref{sec:Euler} we first prove some results for
  the Euler equations, improving the conditions on the time-variable and extending the
  range of allowed H\"older exponent, cf.~Bardos and Titi~\cite{BT2018}. The next result
  in Section~\ref{sec:NSE} is about the viscous case, where we again give a new condition
  in the time variable, beside working in spaces with the same scaling of the known
  results by Cheskidov and Luo~\cite{CL2020}. Finally, in the last and more original
  Section~\ref{sec:half} we treat the viscous problem supplemented with Dirichlet
  conditions in the half-space, proving energy conservation for velocities which are only
  partially H\"older continuous.

% The last result concerns the limiting behavior as the viscosity
% $\nu\to0$ and it is in the same spirit of recent results by Drivas
% and Eyink~\cite{DE2019}. We show how if certain H\"older norms are
% bounded independently of the viscosity, then the solutions of the
% NSE converge to a weak solution of the Euler equations conserving
% the energy.
%
\section{Notation and preliminaries.}
In the sequel we will use  Lebesgue
$(L^{p}(\tore),\|\,.\,\|_{p})$ and Sobolev
$(W^{1,p}(\tore),\|\,.\,\|_{1,p})$ spaces; for simplicity we denote by
$(\,.\,,\,.\,)$ and $\|\,.\,\|$ the $L^{2}$ scalar product and norm
(while the other norms are explicitly indicated).  Moreover, we will
use the Banach space of (uniformly) H\"older continuous functions
$C^{\alpha}(\tore)=C^{0,\alpha}(\tore)$, for $0<\alpha\leq1$, with the
norm
\begin{equation*}
  \|\bu\|_{C^{\alpha}}=\max_{\bx\in\overline{\tore}}|\bu(\bx)|+[u]_{\alpha},
\end{equation*}
where
\begin{equation*}
  [u]_{\alpha}:=
  \sup_{\bx\not=\by}\frac{|\bu(\bx)-\bu(\by)|}{|\bx-\by|^{\alpha}}.       
\end{equation*}
We will focus on space-time properties of functions and what will be
relevant in the sequel is the ``homogeneous behavior'', that is the
one of the H\"older semi-norm $[\,.\,]_{\alpha}$, and we denote by
${\dot C}^{\alpha}$ the space of measurable functions such that this
quantity is bounded\footnote{Note that for bounded domains the uniform
  H\"older property implies also boundedness.}.

We say that
\begin{equation*}
u\in L^{\beta}(0,T;{\dot C}^{\alpha}(\tore)),  
\end{equation*}
if there exists
$f_{\alpha}:[0,T]\to\R^{+}$ such that
\begin{equation*}
  \begin{aligned}
  a)\qquad&  | u(t,x)-u(t,y)|\leq
    f_{\alpha}(t)|\bx-\by|^{\alpha},\qquad\forall\,\bx,\by\in
    \tore,\ \text{for a.e. }\,t\in[0,T],
    \\
b)\qquad &    \int_{0}^{T}f_{\alpha}^{\beta} (t)\,\de
    t<\infty,
  \end{aligned}
\end{equation*}
and note that $f_{\alpha}(t)=[u(t)]_{\alpha}$ for almost all
$t\in[0,T]$. This space will be endowed with the following semi-norm
\begin{equation*}
  \|u\|_{L^{\beta}(0,T;{\dot C}^{\alpha}(\tore))}:=
  \Big[\int_{0}^{T}%(\|u(t)\|_{\infty}^{\beta}+
  f_{\alpha}^{\beta} 
  (t)\,\de t\Big]^{1/\beta}.
\end{equation*}
We define properly the notions of weak solution which we will use in
the sequel. We focus now on the space periodic setting and we will
explain needed changes in the last section.  We denote by $H$ and $V$
the closure of smooth, periodic, divergence-free, and with zero
mean-value vector fields in $L^{2}(\tore)$ and $W^{1,2}(\tore)$. As
space for test functions we use
\begin{equation*}
  \mathcal{D}_T:=\Big\{\bphi\in \bC^{\infty}_{0}([0,T[\times\tore):\
  \dive\bphi=0\Big\}.
\end{equation*}
We define the notion of weak solution in the inviscid case 
  \begin{definition}[Weak solutions for the Euler equations]
  \label{def:Euler-weak-solution}
Let be given $\bv_0\in H$. A measurable function
  $\bve:\,(0,T)\times \tore\to \R^3$ is called a  weak
  solution to the Euler equations~\eqref{eq:NSEbvip} if the
  following hold true:  
\begin{equation*}
  \bve\in L^\infty(0,T;H),
\end{equation*}
and if $\bve$ solves the equations in the weak sense:
\begin{equation}
  \label{eq:distributional-solution2E}
  \int_0^\infty\Big[(\bve,\partial_{t}\bphi)+
  ((\bve\otimes\bve),\nabla\bphi)\Big]\,\de t=-(\bvn_0,\bphi(0)),
\end{equation}
for all $\bphi\in \mathcal{D}_T$.
\end{definition}
We also recall the definition of weak solution for the viscous problem.
\begin{definition}
  \label{def:LH-weak-solution}
  (Space-periodic Leray-Hopf weak solution) Let be given $\bv_0\in
  H$. A measurable function $\bvn:\,(0,T)\times \Omega\to \tore$ is
  called a Leray-Hopf weak solution to the space-periodic
  NSE~\eqref{eq:NSEbvip} if $  \bvn\in L^\infty(0,T;H)\cap L^2(0,T;V);
  $ and the following hold true:
  \\
  % \begin{equation*}
%  \bvn\in L^\infty(0,T;H)\cap L^2(0,T;V);
%\end{equation*}
The function  $\bvn$ solves the equations in the weak sense:
\begin{equation}
  \label{eq:distributional-solution2}
  \int_0^\infty\Big[(\bvn,\partial_{t}\bphi)-\nu(\nabla\bvn,\nabla\bphi)-
  ((\bvn\cdot\nabla)\,\bvn,\bphi)\Big]\,\de t=
  % -\int_0^\infty\langle f,\bphi\rangle\,\de t
  -(\bvn_0,\bphi(0)),
\end{equation}
for all $\bphi\in \mathcal{D}_T$;
\\
It holds the (global) energy inequality 
\begin{equation}
  \label{eq:energy_inequality}
  \frac{1}{2} \|\bvn(t)\|^2_{2}+\nu\int_0^t\|\nabla\bvn(\tau)\|^2_{2}\,
  \de \tau\leq\frac{1}{2}\|\bvn_0\|^2_{2}%+\int_0^t\langle f,\bvn\rangle\,\de
                                %s
  ,\qquad\forall\,t\in[0,T];
\end{equation}
\\
The initial datum is strongly attained 
\begin{equation}
  \label{eq:initial_datum}
  \lim_{t\to0^+}\|\bvn(t)-\bvn_0\|=0.
\end{equation}
\end{definition}
It is well-known that for all $v_{0}\in H$ there exists at least a
Leray-Hopf weak solution in any  time interval $(0,T)$.

The energy inequality~\eqref{eq:energy_inequality} can be rewritten as
an equality
\begin{equation}
  \label{eq:energy_inequality_D}
  \frac{1}{2} \|\bvn(t)\|^2+\nu\int_0^t\int_{\tore}\mbox{\large
    $\epsilon$}[\bvn(\tau,x)]\,\dx 
  \de\tau=\frac{1}{2}\|\bvn_0\|^2,%+\int_0^t\langle f,\bvn\rangle\,\de s
  \qquad\forall\,t\in[0,T],
\end{equation}
where the total energy dissipation  rate is given by 
\begin{equation*}
\mbox{\large  $\epsilon$}[\bv]:=\nu|\nabla\bvn|^{2}+D(\bvn),
\end{equation*}
with $D(\bvn)$ a non-negative distribution (Radon measure); if
$D(\bvn) = 0$, then energy dissipation arises entirely from
viscosity. 
\begin{remark}
  Other equivalent weak formulations can be used (especially for
  the NSE, see Temam~\cite[Ch.~III,\S~3]{Tem1977b}) employing the
  time-derivative as an object at least in $L^{1}(0,T;V^{*})$. This
  requires to find $\bvn\in L^{\infty}(0,T;H)\cap L^{2}(0,T;V)$ with
  $\partial_{t}\bvn\in L^{1}(0,T;V^{*})$ --which in fact it will turn out to
  be in $L^{4/3}(0,T;V^{*})$-- such that the energy inequality is
  satisfied, $\bvn(0)=\bvn_{0}$, and 
\begin{equation*}
%  \label{eq:distributional-solution2}
  \int_0^T\Big[\langle\partial_{t}\bvn,\bphi\rangle+\nu(\nabla\bvn,\nabla\bphi)+
  ((\bvn\cdot\nabla)\,\bvn,\bphi)\Big]\,\de t=0,
  % -\int_0^\infty\langle f,\bphi\rangle\,\de t
%  -(\bvn_0,\bphi(0)),
\end{equation*}
for all divergence-free $\bphi\in \bC^{\infty}_{0}(]0,T[\times\tore)$.
Here $\langle\cdot,\cdot\rangle$ denotes the duality pairing between
$V$ and its dual $V^{*}$. The latter formulation will be particularly
useful when using as test function regularized velocities.
\end{remark}

As for what concerns the Euler equations, we will essentially apply
the same classical strategy as in~\cite{CET1994,DR2000,Eyi1994}),
based on mollifications to justify the calculations.  
%he improved results come from with different
%estimates, which are adapted to Leray-Hopf weak solutions and in the
%simplified setting of classical H\"older continuous functions, instead
%of other function setting.
Note that, \textit{very roughly speaking} the introduction of the
H\"older space $C^{1/3}$ (or  Besov generalizations) is based on
the idea of ``equally distributing'' in the term
\begin{equation*}
  \int_{0}^{T}\int_{\tore}(\bv\cdot\nabla)\,\bv\cdot\bv\dx\de t,
\end{equation*}
a single space derivative, into $1/3$-spatial derivative on all
terms. Despite the approach being the same, results can be greatly
improved for the NSE, since one can take advantage of the already
known $L^{2}$-integrability of the gradient and a different
combination of estimates (which are $\nu$-dependent) could be used.
\section{Mollification and H\"older spaces}
To fully use the features of the H\"older continuity extra-assumption
on the solution, the approximation should be done by a mollification
argument. To this end we fix a symmetric
$\rho\in C^{\infty}_{0}(\R^{3})$ such that
\begin{equation*}
  \rho\geq0,\qquad\text{supp }\rho\subset
  B(0,1)\subset\R^3,\qquad\int_{\R^{3}}\rho(x)\dx=1,  
\end{equation*}
and we define, for $\vep\in(0,1]$, the family of Friedrichs mollifiers
$\rho_{\vep}(x):=\vep^{-3}\rho(\vep^{-1}x)$. Then, for any function
$f\in L^{1}_{loc}(\R^{3})$ we set, by using the usual convolution
(and not the circular convolution made by periodic kernels),
\begin{equation*}
  f_{\vep}(x):=\int_{\R^{3}}\rho_{\vep}(x-y)f(y)\dy=\int_{\R^{3}}\rho_{\vep}(y)f(x-y)\dy. 
\end{equation*}
It turns out that the last integral is evaluated only for
$\{y: |y|<\vep\}\subset]-\pi,\pi[^{3}$, for small $\vep>0$. So if
needed we can restate the definition as
\begin{equation*}
  f_{\vep}(x):=\int_{\tore}\rho_{\vep}(y)f(x-y)\dy. 
\end{equation*}
If 
$f\in L^{1}(\tore)$, then $f\in L^{1}_{loc}(\R^{3})$,  and it turns
out that $f_{\vep}$ is $2\pi$-periodic along the
$x_{j}$-direction, for $j=1,2,3$ since 
\begin{equation*}
  f_{\vep}(x+2\pi e_{j})= \int_{\tore}\rho_{\vep}(y)f(x+2\pi
  e_{j}-y)\dy=\int_{\tore}\rho_{\vep}(y)f(x-y)\dy, 
\end{equation*}
by the periodicity of $f$, where $e_{j}$, $j=1,2,3$, is the unit
vector in the $x_{j}$ direction, and it turns our that
$f_{\vep}\in C^{\infty}(\tore)$.

We report now the basic calculus estimates which we will use in the sequel.
\begin{lemma}
  \label{lem:convolution-Holder}
  Let $\rho$ as above and let $\bu\in {\dot C}^{\alpha}(\tore)\cap
  L^{1}_{loc}(\tore)$, then it follows  
  %\begin{equation*}
    \begin{align}
      \label{eq:conv1}      &   \max_{x\in\tore}
                              |\bu(\bx+\by)-\bu(\bx)|\leq
                              [u]_{\alpha}|\by|^{\alpha},
      \\
      \label{eq:conv2}      &   \max_{x\in\tore} |\bu(x)-\bu_{\vep}(x)|\leq
                              [u]_{\alpha}\, \vep^{\alpha}, 
      \\
      \label{eq:conv3}      &\max_{x\in\tore}|\nabla\bu_{\vep}(x)|\leq C
                              [u]_{\alpha}\,\vep^{\alpha-1}, 
    \end{align}
  %\end{equation*}
    where $C:=\int_{\R^{3}}|\nabla\rho(x)|\dx$.
\end{lemma}
\begin{proof}
  The first one is just the statement of H\"older
  continuity. Concerning the second one we can write that
\begin{equation*}
  \begin{aligned}
    |\bu(\bx)-\bu_{\vep}(\bx)|&=\left|\bu(\bx)-
      \int_{B(0,\vep)}\rho_{\vep}(\by)\bu(\bx-\by)\,\de \by\right|
    \\
    &
    =\left|\int_{B(0,\vep)}\rho_{\vep}(\by)(\bu(\bx)-\bu(\bx-\by))\,\de
      \by\right|
    \\
    &\leq
    [u]_{\alpha}\int_{B(0,\vep)}\rho_{\vep}(\by)|\by|^{\alpha}\,\de
    \by
    \\
    &\leq C \vep^{\alpha}\int_{B(0,\vep)}\rho_{\vep}(\by)\,\de \by,
    \end{aligned}
\end{equation*}
hence the thesis, where we used that $\rho\geq0$.

The third estimate follows by observing that
\begin{equation*}
  \begin{aligned}
    \frac{\partial u_{\vep}(x)}{\partial x_{i}} &= \frac{\partial
    }{\partial x_{i}}
    \frac{1}{\vep^{3}}\int_{\R^{3}}\rho\left(\frac{x-y}{\vep}\right)
    u(y)\dy= \int_{\R^{3}} \frac{1}{\vep^{4}}\frac{\partial
      \rho}{\partial x_{i}}_{\big|\frac{x-y}{\vep}} u(y)\dy
    \\
    &= \frac{1}{\vep}\int_{\R^{3}}g^{i}_{\vep}(x-y)u(y)\dy,
  \end{aligned}
\end{equation*}
where $g^{i}_{\vep}(x):=\frac{1}{\vep^{3}}\frac{\partial
  \rho}{\partial x_{i}}\big(\frac{x}{\vep}\big)$. Note that
$\int_{\R^{3}} g^{i}_{\vep}(x)\,dx=0$, hence 
we can write 
\begin{equation*}
  \begin{aligned}
    \frac{\partial u_{\vep}(x)}{\partial x_{i}} &=\frac{1}{\vep}
    \int_{\R^{3}}g^{i}_{\vep}(x-y)u(y)\dy-\frac{u(x)}{\vep}\int_{\R^{3}}g^{i}_{\vep}(x-y)\dy, 
    \\
    &=\frac{1}{\vep}
    \int_{\R^{3}}\Big(g^{i}_{\vep}(x-y)u(y)-g^{i}_{\vep}(x-y)u(x)\Big)\dy,
  \end{aligned}
\end{equation*}
and then 
\begin{equation*}
  \begin{aligned}
    \left| \frac{\partial }{\partial x_{i}} u_{\vep}(x)\right|&= \frac{1}{\vep}
    \int_{|x-y|<\vep}|g^{i}_{\vep}(x-y)||u(y)-u(x)|\dy,
    \\
    &\leq \frac{[u]_{\alpha}}{\vep}
    \int_{|x-y|<\vep}|g^{i}_{\vep}(x-y)||y-x|^{\alpha}\dy
    \\
    &    \leq \frac{[u]_{\alpha}\vep^{\alpha}}{\vep}
    \int_{|x-y|<\vep}|g^{i}_{\vep}(x-y)|\dy,
  \end{aligned}
\end{equation*}
hence the thesis.
\end{proof}
\subsection{A technical extension on the modulus of continuity}
We observe, that some improvements in the results about energy
conservation can be obtained considering the
following spaces
$C^{\alpha}_{\omega}(\tore)\subset C^{\alpha}(\tore) $ defined
through the
norm
\begin{equation*}
  \|\bu\|_{C^{\alpha}_{\omega}}=\max_{\bx\in\overline{\tore}}|\bu(\bx)|
  + [u]_{\omega,\alpha}    
\end{equation*}
where 
\begin{equation*}
  [u]_{\omega,\alpha}:=
  \sup_{\bx\not=\by}\frac{|\bu(\bx)-\bu(\by)|}{\omega(|\bx-\by|)
    |\bx-\by|^{\alpha}}, 
\end{equation*}
with $\omega:\R^{+}\to \R^{+}$ a non-decreasing function such that
$\lim_{s\to0^+}\omega(s)=0$.  This spaces have been already considered
by Duchon  and Robert~\cite{DR2000} and, more recently in~\cite{BY2020}
for $\alpha=1/3$. See also the Besov counterpart
in~\cite{CCFS2008}. With the same approach used in the previous lemma
one can prove also the following extension of the the properties of
convolutions.
\begin{lemma}
  \label{lem:convolution-Holder2}
  Let $\rho$ as above and let
  $\bu\in {\dot C}^{\alpha}_{\omega}(\tore)\cap L^{1}_{loc}(\tore)$,
  then it follows
  %\begin{equation*}
    \begin{align}
\label{eq:conv1omega}      &   \max_{x\in\tore}
                             |\bu(\bx+\by)-\bu(\bx)|\leq
                             [u]_{\omega,\alpha}\,\omega(|y|)\,|\by|^{\alpha},
      \\
      \label{eq:conv2omega}      &   \max_{x\in\tore}
                                   |\bu(x)-\bu_{\vep}(x)|\leq
                                   [u]_{\omega,\alpha} \,\omega(\vep) \vep^{\alpha},
      \\ 
      \label{eq:conv3omega}      &\max_{x\in\tore}|\nabla\bu_{\vep}(x)|\leq C
       [u]_{\omega,\alpha}\, \omega(\vep)\vep^{\alpha-1}.
    \end{align}
  %\end{equation*}
\end{lemma}
\section{Energy conservation for the Euler equations}
The results in this section concern energy conservation for weak
solutions of the Euler equations:
Theorem~\ref{thm:theorem-energy-holder-Euler} gives a small
improvement in the time integrability, and extends the results to the
full range of $\alpha\in]1/3,1[$, see Bardos and Titi~\cite{BT2018}
and also~\cite{CET1994}. Note that the result about non-conservation of
energy below the exponent $1/3$ in Isett~\cite{Ise2018} are proved with
$L^{\infty}$-bounds in time. Even if the original conjecture
detected the value $\alpha=1/3$ as borderline, we identify
$L^{1/\alpha}(0,T;{\dot C}^{\alpha})$ as a critical space (for all
$\alpha>1/3$).

On the other hand the second result, namely
Theorem~\ref{thm:theorem-energy-holder-Euler-log}, was already known
for $\alpha=1/3$ and here we prove just an extension to the full range
of exponents, see~\cite{BY2020}.
\label{sec:Euler}
 \begin{theorem}
   \label{thm:theorem-energy-holder-Euler}
   Let $\bve$ be a weak solution to the Euler
   equation~\eqref{eq:Ebvip} such that 
\begin{equation}
  \label{eq:2}
  \bve\in   L^{\frac{1}{\alpha}+\delta}(0,T;{\dot
    C}^{\alpha}(\tore)),\quad\text{with } \alpha\in]\frac{1}{3},1[
\end{equation}
and with $\delta>0$. Then, the weak solutions $\bve$ conserves the
energy.
\end{theorem}
\begin{proof}
  We define $\bvee =\rho_{\vep}*\bve$ and consider the following
  equality obtained by testing with
  $\rho_{\vep}*(\rho_{\vep}*\bve)=\rho_{\vep}*\bvee$ (which is
  justified for instance in \cite{BT2018})
  \begin{equation}
    \label{eq:energy-epsilon-Euler}
    \frac{1}{2} \|\bve_{\vep}(t)\|^{2}-
    \frac{1}{2}\|\bve_{\vep}(0)\|^{2}=\int_{0}^{t}\int_{\tore}
    (\bve\otimes\bve)_{\vep}:\nabla\bve_{\vep}\dx \de\tau.
  \end{equation}
  The key observation to handle the convective term from the
  right-hand side is the following (Constantin-E-Titi) commutator-type
  identity, valid for any $u\in L^{2}(\tore)$
\begin{equation}
  \label{eq:commutator-Euler}
  (u\otimes u)_{\vep}=u_{\vep}\otimes u_{\vep}
  +r_{\vep}(u,u)-(u-u_{\vep})\otimes(u-u_{\vep}) ,  
\end{equation}
with
\begin{equation*}
  \br_{\vep}(u,u):=\int_{\tore}\rho_{\vep}(\by)(\delta_{\by}u(\bx)
  \otimes\delta_{\by}u(\bx))\dy, 
\end{equation*}
where (using a now common notation borrowed from~\cite{CET1994}) we 
set
\begin{equation*}
  \delta_{\by}u(\bx):=u(\bx-\by)-u(\bx).
\end{equation*}
We apply the decomposition~\eqref{eq:commutator-Euler} to a mollified
weak solution $\bve_{\vep}$. By integration by parts, it follows that
\begin{equation*}
  \int_{\tore}(\bve_{\vep}\otimes\bve_{\vep}):\nabla\bve_{\vep}\dx=0,
\end{equation*}
so we are reduced to study only the contribution of the remaining two
terms in~\eqref{eq:commutator-Euler}, which are coming from the
decomposition of the right-hand side of
\eqref{eq:energy-epsilon-Euler}.

We start from the last one and we split the integrand as follows
\begin{equation*}
  |\bve-\bvee|^{2}|\nabla\bvee|=
  |\bve-\bvee|^{\eta}|\bve-\bvee|^{2-\eta}|\nabla\bvee|,
  \qquad\text{for  some }0\leq\eta\leq2.
  \end{equation*}
Hence, by using~\eqref{eq:2}, we get  
\begin{equation*}
  \int_{\tore}    |\bve-\bvee|^{2}|\nabla\bvee|\dx\leq
  f_{\alpha}^{1+\eta}(t)\vep^{\alpha\eta+\alpha-1
  }\int_{\tore}|\bve-\bvee|^{2-\eta}    \dx,  
\end{equation*}
and, by H\"older inequality (since the measure of the domain is
finite) we get
$\int_{\tore} |g|^{2-2\eta}\dx\leq C\|g\|_{2}^{2-2\eta}$. Using the
fact that $\|\bvee(t)\|\leq\|\bve(t)\|\leq C$, by the properties of
convolution and the Definition~\ref{def:Euler-weak-solution} of weak
solution, one obtains
\begin{equation*}
\int_{0}^{T}  \int_{\tore}    |\bve-\bvee|^{2}|\nabla\bvee|\dx\de t\leq
C\,\vep^{\alpha\eta+\alpha-1 } \int_{0}^{T}
f_{\alpha}^{1+\eta}(t)\,\de t. 
\end{equation*}
If we want that this term vanishes as $\vep\to0$ we have to satisfy  the
following two facts:
\\
i) to fix $\eta\in[0,2]$ such that $\alpha\,\eta+\alpha-1>0$. The
constraints to be satisfied are
\begin{equation*}
\frac{1-\alpha}{\alpha}<  \eta\leq 2,
\end{equation*}
and such a choice of $\eta$ is possible only if $\alpha>1/3$;
\\
ii) the integral from the right-hand side should be finite, which
means $f_{\alpha}\in L^{1+\eta}(0,T)$ and consequently the inequality
$1+\eta\leq1/\alpha+\delta$ should be verified.

Collecting the two estimates we find that we have to choose $\eta$
such that $1/\alpha-1<\eta\leq1/\alpha-1+\delta$.

The other term arising from the commutator is estimated as follows, by
fixing $\eta$ in the same way as above:
\begin{equation*}
  \begin{aligned}
    |\br_{\vep}(\bve,\bve)|&\leq
    \int_{\R^{3}}\rho_{\vep}(\by)|\bve(x-y)-\bve(x)|^{2}\dy
    \\
    &\leq
    \int_{B(0,\epsilon)}\rho_{\vep}(\by)|\bve(x-y)-\bve(x)|^{\eta}
    |\bve(x-y)-\bve(x)|^{2-\eta}\dy 
    \\
    &\leq f_{\alpha}^{\eta}(t)
    \int_{B(0,\vep)}\rho_{\vep}(\by)|y|^{\alpha
      \eta}|\bve(x-y)-\bve(x)|^{2-\eta}\dy, 
    \\
 &   \leq C(\eta) f_{\alpha}(t)^{\eta} \vep^{\alpha\eta}
    \int_{\tore}\rho_{\vep}(\by)(|\bve(x-y)|^{2-\eta}+|\bve(x)|^{2-\eta})\dy. 
  \end{aligned}
\end{equation*}
Next, 
\begin{equation*}
  \begin{aligned}
    \int_{\tore} \int_{\tore}\rho_{\vep}(\by) \dy 
    |\bve_{\vep}(x)|^{2-\eta}\dx&=\int_{\tore}
    |\bve_{\vep}(x)|^{2-\eta}\dx\leq C(\eta),
  \end{aligned}
\end{equation*}
again by using the H\"older inequality, the fact that $\rho\geq0$ is a
mollifier. While by the properties of the convolution in $L^{p}$ we
have
\begin{equation*}
  \bigg\|
  \int_{\tore}\rho_{\vep}(\by)|\bve(x-y)|^{2-\eta}\dy
  \bigg\|_{L^{1}(\tore)}
  \leq \|
  \rho_{\vep}\|_{1}\||\bve|^{2-\eta}\|_{1}\leq C(\eta),
\end{equation*}
where we used H\"older inequality and the uniform bound in
$L^{2}$. Consequently
\begin{equation*}
  \begin{aligned}
    \left|\int_{\tore}\br_{\vep}(\bve,\bve):\nabla\bve_{\vep}\dx\right|&
    \leq \int_{\tore}|\br_{\vep}(\bve,\bve)| |\nabla\bve_{\vep}|\dx
    \\
    &\leq C\, f_{\alpha}(t)^{1+\eta} \vep^{\alpha \eta+\alpha-1},
\end{aligned}
\end{equation*}
showing that also this term can be treated as the previous one. Hence
we obtain that
\begin{equation*}
  \left|\int_{0}^{t}\int_{\tore}
    (\bve\otimes\bve)_{\vep}:\nabla\bve_{\vep}\dx \de\tau\right|\leq
  C\,\vep^{\gamma } \int_{0}^{T} f_{\alpha}^{1/\alpha+\delta}(t)\,\de t\to0, 
\end{equation*}
under the hypothesis~\eqref{eq:2}, being the integral finite and
$\gamma=\alpha\eta+\alpha-1>0$.
\end{proof}
\begin{remark}
  The result proved means that, apart some ``$\delta$-technical
  conditions'', the energy conservation holds for weak solutions
  $\bve\in L^{1/\alpha}(0,T;{\dot C}^{\alpha})$, for $\alpha>1/3$.
  Observe also that if $\alpha>1/3$, then $1/\alpha<3$, and beside
  having a full-range of scaled space/time results, we have a
  technical-improvement in the time variable with respect to the
  previous results: we can have  $\bve \in L^{q}(0,T;{\dot
    C}^{\alpha}(\tore))$ for some $q<3$, even 
  for $\alpha>1/3$ arbitrarily close to $1/3$.
\end{remark}
By using the slightly smaller space $C^{\alpha}_{\omega}(\tore)$ we can get a result with
a little gain in the time variable. The space $C^{\alpha}_{\omega}(\tore)$ can be
considered as a counterpart of the space $B^{\alpha}_{p ,c(\N)}$, introduced
in~\cite{CCFS2008}, and nevertheless in in the same spirit of the results
in~\cite{BY2019,DR2000} which were limited to the case $\alpha=1/3$. The result which
holds is the following
\begin{theorem}
  \label{thm:theorem-energy-holder-Euler-log}
  Let $\bve$ be a weak solution to the Euler equation~\eqref{eq:Ebvip}
  such that 
\begin{equation}
  \label{eq:3}
  \bve\in 
  L^{1/\alpha}(0,T;{\dot C}^{\alpha}_{\omega}(\tore)), \quad\text{with
}    \alpha\in[\frac{1}{3},1].
\end{equation}
Then $\bve$ conserves the energy.
\end{theorem}
\begin{proof}
  The proof is exactly the same as before, the only difference is the
  following estimation (and an analogous one for the other term in the
  commutator)
  \begin{equation*}
      \int_{0}^{T}  \int_{\tore}    |\bve-\bvee|^{2}|\nabla\bvee|\dx\de
      t\leq \omega(\varepsilon)^{1+\eta}\vep^{\alpha\eta+\alpha-1
      }\int_{0}^{T}      f_{\alpha}^{\eta+1}(t)\de t,
      % \!\int_{\tore}|\bve-\bvee|^{2-\eta}    \dx,   
\end{equation*}
which allows us in this case to fix
$\eta=\frac{1-\alpha}{\alpha}\leq2$, a choice which is possible now
also in the case $\alpha=1/3$. Due to the presence of the term
$\omega(\vep)$, the space-time integral vanishes, as
$\varepsilon\to0$, if $f_{\alpha}\in L^{1/\alpha}(0,T)$.
\end{proof}
\begin{remark}
  Results in this section can be extended to the boundary value
  problem for the Euler equations, with the assumption that
  $\bve\cdot n=0$ at the boundary $\partial\Omega$, with $n$ the unit
  outward vector. This is possible but with a certain amount of
  technicalities, by following word-by-word the approach
  in~\cite{BT2018}. Mainly, the localization does not preserve the
  divergence-free condition and consequently an estimate on the
  pressure should be obtained by solving in the H\"older spaces the
  steady Poisson-Neumann problem
  \begin{equation*}
    \begin{aligned}
      -\Delta p^{E}&=\dive\dive(\bve\otimes\bve)\qquad&&\text{in
      }\Omega, 
      \\
      -n\cdot \nabla p^{E}&=(\bve\otimes\bve) :\nabla
      n\qquad&&\text{on }\partial\Omega, 
      \\
    \end{aligned}
  \end{equation*}
  where the time $t\in[0,T]$ is just a parameter.
\end{remark}
\section{Energy conservation for the Navier-Stokes equations}
\label{sec:NSE}
In this section we focus on the viscous case and recall that the
classical results on scale invariant spaces by
Lady\v{z}henskaya-Prodi-Serrin imply regularity within the class
$\bvn\in L^{r}(0,T;L^{s})$ with $\frac{2}{r}+\frac{3}{s}=1$,
$s>3$. Even the limiting case $s=3$ is a regularity class, but in our
setting it is relevant the other limiting case
\begin{equation*}
  \bvn\in L^{2}(0,T;L^{\infty}(\Omega)).
\end{equation*}
In the light of this observation, full classical regularity (and consequently also energy
conservation) holds true in the class of weak solutions such that $\bvn\in
L^{2}(0,T;C^{0,\alpha}(\overline{\Omega}))$, for all $\alpha\geq0$. In the bounded domain
this implies that conditions involving an exponent for the time-integrability strictly
smaller than $2$ will make the results proved here non-trivial (The role of this limiting
space is also stressed in \cite[Thm.~1.6]{CL2020b} in connection with non-uniqueness).

The main result we prove (improving the known ones in terms of
properties in the time variable) is the following.
\begin{theorem}
\label{thm:theorem-energy-holder}
  Let $\bvn$ be a Leray-Hopf weak solution in $(0,T)\times\tore$ such
  that for all $\lambda\in]0,T[$
  \begin{equation}
    \label{eq:Holder1-3}
    \bvn\in L^{\frac{2}{1+\alpha}}(\lambda,T;{\dot
      C}^{\alpha}(\tore)),\qquad\text{with }\alpha\in]0,1[.
  \end{equation}
Then, $\bvn$ satisfies the energy equality in $[0,T]$.
\end{theorem}
The condition~\eqref{eq:Holder1-3} means that  that there exists
$f:[0,T]\to\R^{+}$ such that 
$| \bvn(t,x)-\bvn(t,y)|\leq
    f_{\alpha}(t)|\bx-\by|^{\alpha}$ for almost all $t\in]0,T]$
    and $f_{\alpha}\in   L^{\frac{2}{1+\alpha}}(\lambda,T)$ 
      for all $\lambda\in]0,T[$,
but it is not excluded that
\begin{equation*}
\int_{0}^{T}f_{\alpha}^{\frac{2}{1+\alpha}}(t)\,\de    t=+\infty,
  \end{equation*}
  which can be restated as  $f_{\alpha}\in L ^{\frac{2}{1+\alpha}}_{loc}(]0,T])$.
\begin{remark}
  Since for all $0<\alpha<1$ we have
  ${\frac{2}{1+\alpha}}<\frac{1}{\alpha}$, the conditions proved in
  Theorem~\ref{thm:theorem-energy-holder} are less restrictive than
  those proved for the Euler equations in the previous section. This
  is due to the fact we have at disposal some additional regularity,
  being $\bvn$ a Leray-Hopf weak solution.

  Note also that the conditions in
  Theorem~\ref{thm:theorem-energy-holder} have the same scaling of the
  recent results from~\cite{CL2020} (recalled
  in~\eqref{eq:Cheskidov-Luo}) in the case $p=\infty$. Nevertheless in
  this case the two results are not completely comparable, since the
  space $L^{\beta}_{w}(0,T)$ is larger than $L^{\beta}(0,T)$, on the
  other hand our results allows for less regularity near $t=0$, using
  the properties of Leray-Hopf solutions, as exploited also in
  Maremonti~\cite{Mar2018}.
\end{remark}
The two limiting cases $\alpha=0$ and $\alpha=1$ in
\eqref{eq:Holder1-3} correspond to
\begin{equation*}
  \begin{aligned}
    &a)\qquad    L^{2}(\lambda,T;L^{\infty}(\Omega)),
    \\
    &b)\qquad     L^{1}(\lambda,T;W^{1,\infty}(\Omega)),
  \end{aligned}
\end{equation*}
and if $\lambda=0$ the case $a)$ is the case covered by the scaling
invariant condition of Lady\v{z}henskaya-Prodi-Serrin~\cite{Gal2000a},
while $b)$ corresponds to the Beale-Kato-Majda criterion. In terms of
scaling (using DeVore diagrams to ``measure'' the regularity of
Sobolev and H\"older spaces) we can associate to the space
$C^{0,\alpha}$ the Sobolev space
$W^{1,\frac{3}{1-\alpha}}\subset C^{0,\alpha}(\overline{\Omega})$ by
Morrey theorem. Results from Theorem~\ref{thm:theorem-energy-holder}
are comparable with
$\nabla v \in
L^{\frac{2}{1+\alpha}}(\lambda,T;L^{\frac{3}{1-\alpha}}(\Omega))$
which is then, in terms of scaling such that
\begin{equation*}
  \frac{2}{\frac{2}{1+\alpha}}+\frac{3}{{\frac{3}{1-\alpha}}}=2,
\end{equation*}
hence with the same scaling-invariant class of regularity
from~\cite{Bei1995a,Ber2002a}.

Moreover, note that the relation between Lady\v{z}henskaya-Prodi-Serrin results and
non-uniqueness has been recently considered in \cite{CL2020b} obtaining as by product
non-uniqueness in the class $L^{3/2-\delta}(0,T;C^{1/3}(\tore)\cap
L^1(0,T;C^{1-\delta'}(\tore))$ for any small $\delta,\delta'>0$. Note that this has
exactly the same scaling as our results for $\alpha=1/3$ and $\alpha=1$, hence showing the
critical role of the condition (at least for very weak solutions of the NSE).
\begin{proof}[Proof of Theorem~\ref{thm:theorem-energy-holder}]
  First, we restrict to the interval $[s,t]$, with $0<s<t\leq T$. We
  use as test function a double smoothed velocity
  $\rho_{\vep}*\bvn_{\vep}:=\rho_{\vep}*(\rho_{\vep}*\bvn)$, obtaining the
  following equality
  \begin{equation}
    \label{eq:energy-epsilon}
    \frac{1}{2} \|\bvn_{\vep}(t)\|^{2}-
    \frac{1}{2}\|\bvn_{\vep}(s)\|^{2}
    +\nu\int_{s}^{t}\|\nabla\bvn_{\vep}(\tau)\|^{2}\,\de 
    \tau=\int_{s}^{t}\int_{\tore}
    (\bvn\otimes\bvn)_{\vep}:\nabla\bvn_{\vep}\dx \de\tau.
  \end{equation}
  To justify the calculations, observe that for weak solution it holds
  that $ \partial_{t}\bv\in L^{4/3}(0,T;V')$, hence the coupling
  \begin{equation*}
    \langle  \partial_{t}\bvn,\rho_{\vep}*(\rho_{\vep}*\bvn)\rangle_{
      L^{4/3}(0,T;V'),L^{4}(0,T;V)},
\end{equation*}
is well defined, being the right-hand side at least in $L^{4}(s,T;V)$.
By means of a further time-mollification one can immediately deduce
that
\begin{equation*}
\int_{s}^{t}\langle
\partial_{t}\bvn,\rho_{\vep}*(\rho_{\vep}*\bvn)\rangle\,\de
\tau=\frac{1}{2}\|\bvn_{\vep}(t)\|^{2}-
\frac{1}{2}\|\bvn_{\vep}(s)\|^{2}.
\end{equation*}
Hence, see for instance Galdi~\cite{Gal2000a},
\begin{equation*}
\lim_{\vep\to0}\int_{s}^{t}\langle
\partial_{t}\bvn,\rho_{\vep}*(\rho_{\vep}*\bvn)\rangle\,\de  
\tau=\frac{1}{2}\|\bvn(t)\|^{2}-\frac{1}{2}\|\bvn(s)\|^{2}.
\end{equation*}
Next, the same classical arguments show also that
\begin{equation*}
  \int_{s}^{t}\int_{\tore}\nabla \bvn_{\vep}: \nabla \bvn_{\vep}\dx
  \de \tau
  \overset{\vep\to0}{\to}
  \int_{s}^{t}\int_{\tore}|\nabla \bvn|^{2}\dx \de \tau.
\end{equation*}
%The key observation to handle the convective term is the following (Constantin-E-Titi) commutator-type identity
%\begin{equation}
%  \label{eq:commutator}
%  (\bv\otimes\bv)_{\vep}=\bvn_{\vep}\otimes\bvn_{\vep}+r_{\vep}(\bv,\bv)
%  -(\bv-\bvn_{\vep})\otimes(\bv-\bvn_{\vep}) , 
%\end{equation}
%with
%\begin{equation*}
%  \br_{\vep}(\bv,\bv):=\int_{\tore}\rho_{\vep}(\by)(\delta_{\by}\bv(\bx)
%\otimes\delta_{\by}\bv(\bx))\dy, 
%\end{equation*}
%where (using a now common notation borrowed from~\cite{CET1994}) we set
%\begin{equation*}
%  \delta_{\by}\bv(\bx):=\bv(\bx-\by)-\bv(\bx).
%\end{equation*}
Next, since by integration by parts
\begin{equation*}
  \int_{\tore}(\bvn_{\vep}\otimes\bvn_{\vep}):\nabla\bvn_{\vep}\dx=0,
\end{equation*}
we are reduced (as in the previous section) to study only the
contribution of the remaining two terms
from~\eqref{eq:commutator-Euler} in the decomposition of the
right-hand side of \eqref{eq:energy-epsilon}. By using the properties
\eqref{eq:conv2}-\eqref{eq:conv3} of the convolution we get
\begin{equation*}
  \begin{aligned}
    & \left|
      \int_{\tore}(\bvn-\bvn_{\vep})\otimes(\bvn-\bvn_{\vep}):\nabla
      \bvn_{\vep}
      \dx\right|\leq\int_{\tore}|\bvn-\bvn_{\vep}|^{2}|\nabla
    \bvn_{\vep}|\dx
    \\
    &\qquad \leq
    \int_{\tore}|\bvn-\bvn_{\vep}|^{1+\alpha}|\bvn-\bvn_{\vep}|^{1-\alpha}|\nabla
    \bvn_{\vep}|^{\alpha}|\nabla \bvn_{\vep}|^{1-\alpha}\dx
    \\
    &\qquad \leq \int_{\tore}|\bvn-\bvn_{\vep}|^{1+\alpha}
    f_{\alpha}(t)^{1-\alpha}\vep^{\alpha(1-\alpha)}
    f_{\alpha}(t)^{\alpha}\vep^{\alpha(\alpha-1)}
    |\nabla\bvn_{\vep}|^{1-\alpha}\dx
    \\
    &\qquad \leq f_{\alpha}(t) \int_{\tore}|\bvn-\bvn_{\vep}|^{1+\alpha}
    |\nabla\bvn_{\vep}|^{1-\alpha}\dx.
\end{aligned}
\end{equation*}
By using the H\"older inequality with exponents $r=\frac{2}{1+\alpha}$ and
$r'=\frac{2}{1-\alpha}$ we get 
\begin{equation*}
  \begin{aligned}
    \left| \int_{\tore}(\bvn-\bvn_{\vep})\otimes(\bvn-\bvn_{\vep}):\nabla
      \bvn_{\vep} \dx\right|&    \leq
    f_{\alpha}(t)\|\bvn-\bvn_{\vep}\|^{1+\alpha}\|\nabla\bvn_{\vep}\|^{1-\alpha}
    \\
    &      \leq
   (2\|v_{0}\|)^{1+\alpha}  f_{\alpha}(t)\|\nabla\bvn_{\vep}\|^{1-\alpha},
  \end{aligned}
\end{equation*}
where in the last step we used elementary properties of convolutions
and the energy inequality to get
$\|\bvn_{\vep}(t)\|\leq\|\bvn(t)\|\leq\|\bv_{0}\|$.

This proves, by using again H\"older inequality, now with respect to
the time variable, that 
\begin{equation*}
  \begin{aligned}
      \left|\int_{s}^{t} \int_{\tore}(\bvn-\bvn_{\vep})\otimes(\bvn-\bvn_{\vep}):\nabla 
      \bvn_{\vep} \dx\de \tau\right|&\leq C \int_{s}^{T}
      f_{\alpha}(\tau)\|\nabla\bvn_{\vep}(\tau)\|_{2}^{1-\alpha}\,\de \tau
      \\
      &\leq C\|f_{\alpha}\|_{\frac{2}{2+\alpha}}\Bigg(\int_{s}^{T}
      \|\nabla\bvn_{\vep}(\tau)\|_{2}^{2}\,\de \tau\Bigg)^{1-\alpha}
      \\
      &\leq C\|f_{\alpha}\|_{\frac{2}{2+\alpha}}\Bigg(\int_{s}^{T}
      \|\nabla\bvn(\tau)\|_{2}^{2}\,\de \tau\Bigg)^{1-\alpha},
    \end{aligned}
  \end{equation*}
  which is uniform in $\vep>0$, since
  $\nabla \bvn\in L^{2}(0,T;L^{2}(\tore))$. Concerning the
  \textit{remainder} term $\br_{\vep}(\bvn,\bvn)$ we can use the same
  argument based on Young theorem on convolutions to prove that
\begin{equation*}
  \begin{aligned}
    \left|\int_{s}^{t}\int_{\tore}\br_{\vep}(\bvn,\bvn):\nabla\bvn_{\vep}\dx\de\tau
    \right|\leq C_{1}
    \|f_{\alpha}\|_{\frac{2}{2+\alpha}}\Bigg(\int_{s}^{T}
    \|\nabla\bvn(\tau)\|_{2}^{2}\,\de \tau\Bigg)^{1-\alpha}.
  \end{aligned}
\end{equation*}
The above calculations prove that the family of functions
$F_{\vep}:(0,T]\to\R^{+}$, indexed by $\vep>0$, and defined by 
\begin{equation*}
  F_{\vep}(\tau):=\int_{\tore}(\bvn(\tau,\bx) \otimes\bvn(\tau,\bx)
  )_{\vep}:\nabla\bvn_{\vep}(\tau,\bx) \dx\qquad \forall\, \tau\in]0,T], 
\end{equation*}
is such that for any fixed $\nu>0$ and  for any $s\in(0,T)$ 
\begin{enumerate}
\item $|F_{\vep}(\tau)|$ is uniformly bounded in $(s,T)$ by the
  function 
  $G(\tau):=C_{2} f_{\alpha}(\tau)\|\nabla
  \bv(\tau)\|_{2}^{1-\alpha}$ which belongs to $L^{1}(s,T)$ by hypothesis;
\item Due to the a.e. in $(0,T)$ convergence of $\bvn_{\vep}(\tau)$
  towards $\bvn(\tau)\in V$ (valid being $v\in L^{2}(0,T;V)$\,) it
  follows that
  \begin{equation*}
    \lim_{\vep\to0}F_{\vep}(\tau)=\int_{\tore}(\bvn(\tau,x)\otimes
    \bvn(\tau,x)):\nabla  \bvn(\tau,x)\dx=0,\quad  \text{for  a.e. }  \tau\in(0,T).
  \end{equation*}
\end{enumerate}
By using the Lebesgue dominated convergence theorem, the two
properties imply that, 
along a sub-sequence,
\begin{equation*}
  \lim_{\vep\to0}  \int_{s}^{t}F_{\vep}(\tau)\,\de \tau=\lim_{\vep\to0}
  \int_{s}^{t}\int_{\tore}(\bvn(\tau,x)\otimes \bvn(\tau,x))_{\vep}:\nabla 
  \bvn_{\vep}(\tau,x)\dx\de \tau=0,
\end{equation*}
for all $t\in]s,T]$, hence finally proving that 
\begin{equation*}
  \frac{1}{2}\|\bvn(t)\|^{2}+\nu\int_{s}^{t}\|\nabla\bvn(\tau)\|^{2}\,\de \tau=
  \frac{1}{2}\|\bvn(s)\|^{2},\qquad\forall\,s,t\text{ s.t. } 0<s< t\leq T.
\end{equation*}
Next, we take a sequence $\{s_{m}\}$ of strictly positive times such
that $s_{m}\to0$. 
By the definition of weak solution it holds that: i) 
\begin{equation*}
\lim_{m\to+\infty}  \|\bv(s_{m})\|=\|\bv_{0}\|,
\end{equation*}
being the initial   datum strongly attained; ii)
\begin{equation*}
\lim_{m\to+\infty}  \int_{s_{m}}^{t}\|\nabla\bvn(\tau)\|^{2}\,\de \tau=
  \int_{0}^{t}\|\nabla\bvn(\tau)\|^{2}\,\de \tau,
\end{equation*}
by the absolute continuity of the time-integral of
$\|\nabla \bv(t)\|^{2}$. Hence, passing to the limit as $m\to+\infty$
in the energy equality over $[s_{m},t]$ is justified and we get that
  \begin{equation*}
 \frac{1}{2}\|\bvn(t)\|^{2}+\nu
 \int_{0}^{t}\|\nabla\bvn(\tau)\|^{2}\,\de \tau= 
 \frac{1}{2}\|\bvn_{0}\|^{2},\qquad \forall\, t\in[0,T].
\end{equation*}
\end{proof}
\section{Energy conservation in presence of boundaries:
  The half-space case} 
\label{sec:half}
In this section we prove a result concerning the energy conservation
in H\"older spaces, in the case of a domain with boundary, and with
homogeneous Dirichlet conditions. The presence of a solid boundary
makes the problem more complex: obtaining smooth and with compact
support approximations is generally subject to a localization process,
which is not preserving the divergence-free constraint and requires
the derivation of appropriate estimates on the pressure. Correction of
the divergence is a non-local process, which seems at present not
easily obtainable preserving all the pointwise estimates, as needed to
take benefit the extra H\"older-continuous assumptions on the
solution. In the case of the NSE this is more technical than for the
Euler equations: in the time-dependent case the estimates in a bounded
domain are generally obtainable through the solution of a linear
evolution Stokes problem
\begin{equation*}
  \begin{aligned}
    \partial_{t}\bvn-\nu\Delta\bvn+\nabla
    p&=-(\bvn\cdot\nabla)\,\bvn\qquad &(t,\bx)\in (0,T)\times \Omega,
    \\
    \dive\bvn&=0\qquad &(t,\bx)\in (0,T)\times \Omega,
    \\
    \bvn&=0\qquad &(t,\bx)\in (0,T)\times \partial\Omega,
    \\
    \bvn(0,\bx)&=\bvn_0(\bx) \qquad &\bx\in \Omega.
  \end{aligned}
\end{equation*}
This requires the construction of a strong solution $(v,p)$, with
right-hand side $-(\bvn\cdot\nabla)\,\bvn$; the pressure cannot be
constructed time-by time as solution of a Poisson problem,
(see~\cite[p.~247]{Soh2001}), as for the Euler equations. Even having
$\bvn\in L^{\beta}(0,T;C^{\alpha}(\Omega))$ (with $\beta<2$ not to
fall into classical regularity classes) will not produce a time-integrable
right-hand side, since $\nabla\bvn\in L^{2}(0,T;L^{2}(\Omega))$.

In addition, other possibilities of constructing directly
divergence-free approximations (as for instance in
Masuda~\cite[Appendix]{Mas1984}) seem at the moment not to produce
interesting results and --up to our knowledge-- the results of energy
conservation in presence of boundaries (in the H\"older/Besov case) is
always subject to extra assumptions on the pressure. In this respect
we wish to mention the results of~\cite{DN2019} where, to study weak
limits, conditions uniform in $\nu\in]0,1]$ are identified. See
also~\cite{BTW2019} with additional conditions on the pressure and the
results in~\cite{DN2018} with an extra $L^{3}$-equi-continuity near
boundary.

In this paper we prove a new result in the half-space, by using a
suitable decomposition of variables, into horizontal and vertical
ones.  We consider the NSE in the half-space
$\R^{3}_{+}:=\left\{x\in \R^{3}: x_{3}>0\right\}$ with vanishing
Dirichlet condition $v=0$ on $\{x\in \R^{3}: x_{3}=0\}$.  In the
sequel denote by $H$ and $V$ the closure of smooth, periodic,
divergence-free and with zero mean-value vector fields in
$L^{2}(\R^{3}_{+})$ and ${\dot W}^{1,2}_{0}(\R^{3}_{+})$. The notion
of Leray-Hopf weak solution remains the same as in
Definition~\ref{def:LH-weak-solution}, just changing the function
spaces and taking test functions with compact support in
$[0,T)\times \R^{3}_{+}$.

The particular geometry of the domain will make the problem tractable,
by means of a natural splitting of variables and unknowns.  To this
end, any $x=(x_{1},x_{2},x_{3})\in \R^{3}$ will be considered as
$x=(x_{h},x_{3})$ and the same decomposition will be applied to the
vector valued functions $u=(u_{h},u_{3})$. We will also denote by
$\nabla_{h}u:=(\partial_{x_{1}}u,\partial_{x_{2}}u) $ and
$\dive_{h}u_{h}:=\partial_{x_{1}}u^{1}+\partial_{x_{2}}u^{2}$ the
gradient and the divergence in the horizontal variables, respectively.
\begin{remark}
  Different from what done in~\cite{RRS2018} for the Euler equations,
  for the NSE the reflection principle does not simply apply. On the
  other hand it is well-known that in the case of \textit{Navier
    boundary conditions in the half -space} the reflection principle
  can be used, giving a direct way to adapt results known in the
  half-space to the full-space.
\end{remark}
We introduce the following class of functions:
\begin{definition}
  We say that $u\in {\dot C}_{\omega,hor}(\R^{3}_{+})$ if for almost
  all $z_{3}>0$
  \begin{equation*}
    \sup_{x_{h}\not=y_{h}}{|u(x_{h},z_{3})-u(y_{h},z_{3})|}\leq{\omega(|x_{h}-y_{h}|)}, 
  \end{equation*}
  with $\omega:\R^{+}\to \R^{+}$  a non-decreasing function such
  that $\lim_{s\to0^+}\omega(s)=0$. This means that, for almost all
  $z_{3}>0$, the function $u$ is uniformly continuous with modulus of
  continuity $\omega(\,.\,)$, with respect to the horizontal variables
  $x_{h}$.
\end{definition}
The main result of this section is the following one:
\begin{theorem}
  Let $\bvn$ be a Leray-Hopf weak solution of the NSE in the half
  space, such that
\begin{equation}
  \label{eq:hyp-half-space}
  \bvn\in L^{2}(0,T;{\dot C}_{\omega,hor}(\R^{3}_{+})).
\end{equation}
Then, $\bvn$ conserves the energy.
\end{theorem}
\begin{remark}
  The result can be interpreted as (especially near to a flat
  boundary) energy is conserved if some regularity of the translations
  parallel to the boundary is assumed. No conditions on the pressure
  or on vertical increments are required.
  
  The same result applies also to a channel flow, that is if the
  domain is $\R^{2}\times(0,1)$, or a periodic strip
  $\T^{2}\times(0,1)$, as studied in Robinson, Rodrigo, and
  Skipper~\cite{RRS2018} with a reflection principle.  Note that the
  regularity assumed in the time variable is more restrictive than in
  the previous section.
\end{remark}
The proof of the main result of this section is based on the use of a
partial mollification. We fix a symmetric
$\widetilde{\rho}\in C^{\infty}_{0}(\R^{2})$ such that
\begin{equation*}
  \widetilde{\rho}\geq0,\qquad\text{supp }\widetilde{\rho}\subset
  B(0,1)\subset
  \R^{2},\qquad\int_{\R^{2}}\widetilde{\rho}(x_{h})\dx_{h}=1,   
\end{equation*}
and we define, for $\vep\in(0,1]$,
$\widetilde{\rho}_{\vep}(x_{h}):=\vep^{-2}\widetilde{\rho} (\vep^{-1}x_{h})$.
\begin{definition}[Horizontal mollification]
  We set, for $u\in L^{p}(K\times\R_{+})$, for all $K\subset\subset\R^{2}$ 
  \begin{equation*}
    \widetilde{u}_{\vep}(x_{h},x_{3}):=
    (\widetilde{\rho}_{\vep}*_{h}u)(x_{h},x_{3})=
    \int_{\R^{2}}\rho_{\vep}(y_{h})u(x_{h}-y_{h},x_{3})\,\de x_{h}, 
  \end{equation*}
  that is a mollification only with respect to the horizontal
  variables $x_{h}$.
\end{definition}

% \marginpar{write in terms of modulus of cont}
\begin{lemma}
  We have the following properties: Let $\widetilde{u}_{\vep}$ be the
  horizontal mollification of
  $u\in W^{1,2}_{0}(\R^{3}_{+})\cap {\dot
    C}_{\omega,hor}(\R^{3}_{+})$, then it holds
  \begin{align}
    &    \widetilde{u}_{\vep}(x_h,0)=0;
    \\
    & \text{if  $\dive u=0$,
      then $\dive \widetilde{u}_{\vep}=0$};
    \\
    &
      \sup_{x\in\R^{3}_{+}} |\bu(x)-\bu_{\vep}(x)|\leq
      \omega(\vep).
 \end{align}
\end{lemma}
\begin{proof}
  The first statement comes from $u(x_{h},0)$, which implies directly
  $\widetilde{\rho}*_{h}u(x_{h},0)=0$, being the integral evaluated in
  a set where the function $u$ vanishes. The divergence of the
  mollified function can be evaluated explicitly, to prove that is
  preserved by the horizontal mollification.

  The third property follows by observing that 
  \begin{equation*}
  \begin{aligned}
    |\bu(\bx)-\widetilde{\bu}_{\vep}(\bx)|&=\left|\bu(\bx_{h},x_{3})-
      \int_{\R^{2}}\widetilde{\rho}_{\vep}(\by_{h})\bu(\bx_{h}-\by_{h},x_{3})\,\de
      \by_{h}\right| 
    \\
    &
    =\left|\int_{B(0,\vep)}\widetilde{\rho}_{\vep}(\by)(\bu(\bx_{h},x_{3})-
      \bu(\bx_{h}-\by_{h},x_{3})\,\de \by_{h}\right|
    \\
    &\leq 
    \int_{B(0,\vep)}\widetilde{\rho}_{\vep}(\by)\,\omega(|\by_{h}|)\,\de
    \by_{h}
    \\
    &\leq \omega(\vep)
\int_{\R^2}\widetilde{\rho}_{\vep}(\by)\,\de
    \by_{h},
  \end{aligned}
\end{equation*}
hence the thesis, where we used that $\widetilde{\rho}\geq0$ and the
total mass equals one.
\end{proof}
The main fact  is that if $\bvn$ is a Leray-Hopf weak solution,
then $\bvnh_{\vep}$ is still vanishing at the boundary,
divergence-free, and it a smooth in the space variables $x_{h}$. To
justify its use as a test functions, a further smoothing in the
time-variable will be needed. As in the previous section we fix
$0<s<t\leq T$, to work (for the moment) in
the fixed time interval
$[s,t]\subset]0,T]$. Next, we define
$u_{\kappa}:\,[s,t]\times \R^{3}_{+}\to\R^{3}$ by
\begin{equation*}
u_{\kappa}(\sigma,x):=\int_{s}^{t}J_{\kappa}(\sigma-\tau)
  u(\tau,x)\,\de \tau\qquad\text{with } 0<\kappa<t-s, 
\end{equation*}
which is obtained through another mollification with respect to the
time-variable (this time by means of the non-negative, even with
respect to $t=s$, smooth and compactly supported function
$J_{\kappa}$, such that $\int_{s-1}^{s+1}J(\tau)\,\de\tau=1$, e.g. one
can use $J(\tau)=\rho(\tau-s)$. )

Then, we will use as test function in the weak formulation of the NSE
the function
\begin{equation*}
  \bvnh_{\vep,\vep,\kappa}(\sigma,x):=\int_{s}^{t}J_{\kappa}(\sigma-\tau)
  \big(\widetilde{\rho}_{\vep}*(\widetilde{\rho}_{\vep}*\bvn)\big)(\tau,x)\,\de 
  \tau.  
\end{equation*} 
Observe that, the use as test function is justified, because it
is zero at the boundary, divergence-free and, from  
\begin{equation*}
  \bvnh_{\vep,\vep,\kappa}\in
  W^{1,\infty}(s,t;V')\cap L^{\infty}(s,t;V), 
\end{equation*}
the following estimates
\begin{equation*}
  \begin{aligned}
&\left|    \int_{s}^{t}<\bvn, \partial_{t}\bvnh_{\vep,\vep,\kappa}>\de
    \tau\right|\leq \|\partial_{t}\bvn\|_{L^{4/3}(V')}\|
  \bvnh_{\vep,\vep,\kappa}\|_{L^{4}(V)}, 
    %\quad\text{and}\quad
  \\
  &\left|\int_{s}^{t}\int_{\R^{3}_{+}}(\bvn\cdot\nabla)\,
    \bvnh_{\vep,\vep,\kappa}\cdot \bvn\dx\de\tau\right|\leq
  C\|\bvn\|_{L^{8/3}(L^{4})}\|\nabla
  \bvnh_{\vep,\vep,\kappa}\|_{L^{\infty}(L^{2})},
\end{aligned}
\end{equation*}
hold true. This shows that the left-hand sides are properly
defined. Moreover, by using the regularity of solutions to the linear
evolution Stokes problem, we can infer (as in Sohr and von
Wahl~\cite{SvW1986}) that the pressure can be selected such that
$\nabla p\in L^{5/4}((s_{n},t)\times \R^{3}_{+})$ for a sequence of
times $s_{n}>0$ such that $s_{n}\to0$. These are chosen such that
$v(s_{n})\in H^{1}_{0}(\R^{3}_{+})$. This shows also that
\begin{equation*}
  \int_{s}^{t}\int_{\R^{3}_{+}}\nabla p  \cdot
  \bvnh_{\vep,\vep,\kappa}\dx\de\tau=0, 
\end{equation*}
since
\begin{equation*}
  \bvnh_{\vep,\vep,\kappa}\in L^{
    \infty}(s,t;L^{2}(\R^{3}_{+})\cap
  L^{\infty}(s,t;{\dot W}^{1,2}_{0}(\R^{3}_{+}))
  \subset L^{\infty}(s,t;L^{2}(\R^{3}_{+})\cap
  L^{6}(\R^{3}_{+})).
\end{equation*}
This finally proves that, by using the so-called Hopf-lemma 
(cf. e.g.~\cite{Gal2000a})  that
\begin{equation*}
  \begin{aligned}
  &  \int_{s}^{t}<\bvnh_{\vep},
    \partial_{t}\bvnh_{\vep,\kappa}>\de\tau-\nu
    \int_{s}^{t}\int_{\R^{3}_{+}} 
    \nabla\bvnh_{\vep}: \nabla\bvnh_{\vep,\kappa}\dx\de\tau+
    \\
    &\qquad \int_{s}^{t}
    \int_{\R^{3}_{+}}\widetilde{(\bvn\otimes\bvn)}_{\vep}:
    \nabla\bvnh_{\vep,\kappa}\dx\de\tau
    =\int_{\R^{3}_{+}}\bvnh_{\vep}(t)\cdot
    \bvnh_{\vep,\kappa}(t)\dx-\int_{\R^{3}_{+}}\bvnh_{\vep}(s)\cdot
    \bvnh_{\vep,\kappa}(s)\dx.
  \end{aligned}
\end{equation*}
Note that 
\begin{equation*}
    \int_{s}^{t}<\bvnh,
    \partial_{t}\bvnh_{\vep,\vep,\kappa}>\de\tau =
\int_{\R^{3}_{+}}    \int_{s}^{t}\int_{s}^{t}J'(\sigma-\tau)\bvnh_{\vep}(\sigma)
    \bvnh_{\vep}(\tau)\,\de\sigma \de\tau\dx=0,
\end{equation*}
being the kernel $J'$ odd with respect to $t=s$. By a standard
limiting procedure $\kappa\to0$, by using the weak
$L^{2}$-continuity of Leray-Hopf weak solutions,  and with the fact that
$\int_{s}^{s+\kappa}J_{\kappa}(\tau)\,\de\tau=\frac{1}{2}$, one
obtains that
\begin{equation*}
  \int_{\R^{3}_{+}}\bvnh_{\vep}(t)\cdot
  \bvnh_{\vep,\kappa}(t)\dx-\int_{\R^{3}_{+}}\bvnh_{\vep}(s)\cdot
  \bvnh_{\vep,\kappa}(s)\dx\overset{\kappa\to0}\to\frac{1}{2}\|\bvnh_{\vep}(t)\|^{2}-\frac{1}{2}\|\bvnh_{\vep}(s)\|^{2},
\end{equation*}
and also trivially
$ \int_{s}^{t}\int_{\R^{3}_{+}} \nabla\bvnh_{\vep}:
\nabla\bvnh_{\vep,\kappa}\dx\de\tau\to
\int_{s}^{t}\|\nabla\bvnh_{\vep,\kappa}(t)\|^{2}\de\tau$, since
$\bvnh_{\vep,\kappa}\to \bvnh_{\vep}$ at least in
$L^{2}((s,t)\times \R^{3}_{+})$.  To handle the nonlinear term observe
that for any Banach space $X$, if $u\in L^{q}(s,t;X)$, then
$u_{\kappa}\to u$ strongly in $L^{q}(s,t;X)$, hence if the limiting
object is finite, then it would follow
\begin{equation}
  \label{eq:limit-nonlinear}
  \int_{s}^{t}
  \int_{\R^{3}_{+}}\widetilde{(\bvn\otimes\bvn)}_{\vep}:
  \nabla\bvnh_{\vep,\kappa}\dx\de\tau\overset{\kappa\to0}{\to}
  \int_{s}^{t} 
  \int_{\R^{3}_{+}}\widetilde{(\bvn\otimes\bvn)}_{\vep}:
  \nabla\bvnh_{\vep}\dx\de\tau.
\end{equation}
To show the validity of the limit~\eqref{eq:limit-nonlinear} we write
as usual
\begin{equation}
  \label{eq:commutator-Euler-h}
  \widetilde{ (u\otimes u)}_{\vep}=\widetilde{u}_{\vep}\otimes
  \widetilde{u}_{\vep}+\widetilde{r}_{\vep}(u,u) 
  -(u-\widetilde{u}_{\vep})\otimes(u-\widetilde{u}_{\vep}) , 
\end{equation}
with
\begin{equation*}
  \widetilde{\br}_{\vep}(u,u):=\int_{\tore}\widetilde{\rho}_{\vep}(\by)
  (\widetilde{\delta}_{\by}u(\bx)  
\otimes\widetilde{\delta}_{\by}u(\bx))\dy, 
\end{equation*}
where we set
\begin{equation*}
  \widetilde{\delta}_{\by}u(\bx):=u(\bx-\by_{h})-u(\bx).
\end{equation*}
We first handle the resulting integral
$ \int_{\tore}(\bvnh_{\vep}\otimes\bvnh_{\vep}):\nabla\bvnh_{\vep}\dx$
which we will prove to be finite by the following observations: for
all $k\geq1$ we have that
$\nabla_{h}\bvnh_{\vep}=(\partial_{x_{1}}\bvnh_{\vep},
\partial_{x_{2}}\bvnh_{\vep})\in L^{2}(0,T;H^{k}(\R^{3}_{+}))$. This implies
\begin{equation*}
  \partial_{x_{3}}\bvnh_{\vep}^{3}=-\dive_{h} \bvnh_{\vep}^{h}\in
  L^{2}(0,T;H^{k}(\R^{3}_{+})),\qquad \text{ for all }k\geq1,
\end{equation*}
hence the third component $\bvnh_{\vep}^{3} $ is smooth in the space
variables.  Observe next that by using  the decomposition
\begin{equation*}
  (\bvnh_{\vep}\cdot\nabla)\bvnh_{\vep}\cdot \bvnh_{\vep}=
  (\bvnh_{\vep}^{h}\cdot\nabla_{h}) \bvnh_{\vep}\cdot \bvnh_{\vep}+
  \bvnh_{\vep}^{3}\, 
  \partial_{3}\bvnh_{\vep}\cdot \bvnh_{\vep},
\end{equation*}
then
\begin{equation*}
  \begin{aligned}
\left|    \int_{s}^{t}\int_{\R^{3}_{+}}
    (\bvnh_{\vep}\cdot\nabla)\bvnh_{\vep}\cdot \bvnh_{\vep}\dx \de
    \tau\right|&=\left|    \int_{s}^{t}\int_{\R^{3}_{+}}
(\bvnh_{\vep}\otimes\bvnh_{\vep})
  :\nabla\bvnh_{\vep}\dx \de\tau\right|
  \\
  &\leq     \int_{s}^{T}
    \|\bvnh_{\vep}^{h}\|_{2}\|\nabla_{h}\bvnh_{
    \vep}\|_{\infty}\|\bvnh_{\vep}\|_{2}+
    \|\bvnh_{\vep}^{3}\|_{\infty}\| \partial_{3}\bvnh_{\vep}\|_{2}\|
    \bvnh_{\vep}\|_{2}\,\de\tau
    \\
    &\leq C(\vep^{-1})\|v\|_{L^{\infty}(H)\cap L^{2}(V)}^{3}.
  \end{aligned}
\end{equation*}
This shows that, for any fixed $\vep>0$,
\begin{equation*}
  \int_{s}^{t}
  \int_{\tore}(\bvnh_{\vep}\otimes\bvnh_{\vep})
  :\nabla\bvnh_{\vep,\kappa}\dx\de
  t\overset{\kappa\to0}{\to}  \int_{s}^{t}
  \int_{\tore}(\bvnh_{\vep}\otimes\bvnh_{\vep}):\nabla\bvnh_{\vep}\dx\de\tau,
\end{equation*}
and, by integration by parts (which is legitimate since the integral
exists), it follows that 
\begin{equation*}
    \int_{s}^{t}
    \int_{\tore}(\bvnh_{\vep}\otimes\bvnh_{\vep}):\nabla\bvnh_{\vep}
    \dx\de\tau=0. 
\end{equation*}
Next, by using the same argument as before, we estimate the first term in the
decomposition as follows
\begin{equation*}
  \begin{aligned}
     \left|
      \int_{\R^{3}_{+}}(\bvn-\bvnh_{\vep})\otimes(\bvn-\bvnh_{\vep}):\nabla
      \bvnh_{\vep}
      \dx\right|
    %\\
    % &\qquad
    &\leq
    \int_{\R^{3}}|\bvn-\bvnh_{\vep}|\,|\bvn-\bvnh_{\vep}|\,|\nabla
    \bvnh_{\vep}|\dx
    \\
    &\leq     f_{\omega}(t)\omega(\vep)
    \|\bvn-\bvnh_{\vep}\|_{2}
   \|\nabla\bvnh_{\vep}\|_{2}
    \\
    &\leq    2 f_{\omega}(t)\omega(\vep)
    \|\bvn\|_{2}
   \|\nabla\bvn\|_{2}.
    % &\qquad \leq f_{\alpha}(t) \int_{\tore}|\bvn-\bvn_{\vep}|^{1+\alpha}
%    |\nabla\bvn_{\vep}|^{1-\alpha}\dx.
\end{aligned}
\end{equation*}
This finally shows that, being the integral finite, then the limit
process in $\kappa\to0$ is valid also in the other term coming from
the commutator decomposition; moreover
\begin{equation*}
  \begin{aligned}
&    \left| \int_{s}^{t}
      \int_{\R^{3}_{+}}(\bvn-\bvnh_{\vep})\otimes(\bvn-\bvnh_{\vep}):\nabla
      \bvnh_{\vep} \dx \de\tau\right|
    \\
    &\qquad \leq\omega(\vep)\|\bvn_{0}\|^{2}\left[\int_{s}^{T}f_{\omega}^{2}(t)\,\de
      t\right]^{1/2} \left[\int_{s}^{T}\|\nabla\bvn (t)\|^{2}\,\de
      t\right]^{1/2}
\\
    &\qquad     \leq C\omega(\vep).
  \end{aligned}
\end{equation*}
The term arising from $\widetilde{r}_{\vep}(u,u)$ is treated in the
same way, hence showing that, for each fixed $\vep>0$, the following
equality holds true 
  \begin{equation*}
    \frac{1}{2}\|\bvnh_{\vep}(t)\|^{2}+\nu\int_{s}^{t}
    \|\nabla\bvnh_{\vep}(\tau)\|^{2} \de \tau-\int_{s}^{t}
    \int_{\R^{3}_{+}}\widetilde{(\bvn\otimes\bvn)}_{\vep}
    :\nabla\bvnh_{\vep}\dx\de \tau=
    \frac{1}{2}\|\bvnh_{\vep}(s)\|^{2}. 
\end{equation*}
Moreover, the bounds proved show also that in the limit $\vep\to0$
\begin{equation*}
\left|    \int_{s}^{t} 
  \int_{\R^{3}_{+}}\widetilde{(\bvn\otimes\bvn)}_{\vep}:
  \nabla\bvnh_{\vep}\dx\de\tau\right|\leq
C\omega(\vep)\overset{\vep\to0}{\to}0, 
\end{equation*}
due to the
hypothesis~\eqref{eq:hyp-half-space}. This proves that
  \begin{equation*}
    \frac{1}{2}\|\bvn(t)\|^{2}+\nu\int_{s}^{t}
    \|\nabla\bvn(\tau)\|^{2} \de \tau=
    \frac{1}{2}\|\bvn(s)\|^{2},\qquad
    \text{for all } 0<s<t\leq T. 
\end{equation*}
The arbitrariness of $s,t$ and the strong
limit~\eqref{eq:initial_datum} is used to end the proof.
%a
\section*{Acknowledgments}

The author acknowledges support by INdAM-GNAMPA.
\def\ocirc#1{\ifmmode\setbox0=\hbox{$#1$}\dimen0=\ht0 \advance\dimen0
  by1pt\rlap{\hbox to\wd0{\hss\raise\dimen0
  \hbox{\hskip.2em$\scriptscriptstyle\circ$}\hss}}#1\else {\accent"17 #1}\fi}
  \def\polhk#1{\setbox0=\hbox{#1}{\ooalign{\hidewidth
  \lower1.5ex\hbox{`}\hidewidth\crcr\unhbox0}}} \def\cprime{$'$}

\end{document}